\pdfoutput=1
\RequirePackage{ifpdf}
\ifpdf % We are running pdfTeX in pdf mode
\documentclass[pdftex]{sigma}
\else
\documentclass{sigma}
\fi

\begin{document}

\allowdisplaybreaks

\renewcommand{\PaperNumber}{055}

\FirstPageHeading

\ShortArticleName{$\mathfrak{spo}(2|2)$-Equivariant Quantizations on the Supercircle $S^{1|2}$}

\ArticleName{$\boldsymbol{\mathfrak{spo}(2|2)}$-Equivariant Quantizations
\\
on the Supercircle $\boldsymbol{S^{1|2}}$}

\Author{Najla MELLOULI~$^\dag$, Aboubacar NIBIRANTIZA~$^\ddag$ and Fabian RADOUX~$^\ddag$}

\AuthorNameForHeading{N.~Mellouli, A.~Nibirantiza and F.~Radoux}

\Address{$^\dag$~University of Sfax, Higher Institute of Biotechnology,
\\
\hphantom{$^\dag$}~Route de la Soukra km 4, B.P.
n\textsuperscript{o} 1175, 3038 Sfax, Tunisia}
\EmailD{\href{mailto:mellouli@math.univ-lyon1.fr}{mellouli@math.univ-lyon1.fr}}

\Address{$^\ddag$~University of Li\`ege, Institute of Mathematics,
\\
\hphantom{$^\ddag$}~Grande Traverse, 12 - B37, B-4000 Li\`ege, Belgium}
\EmailD{\href{mailto:a.nibirantiza@doct.ulg.ac.be}{a.nibirantiza@doct.ulg.ac.be},
\href{mailto:fabian.radoux@ulg.ac.be}{fabian.radoux@ulg.ac.be}}

\ArticleDates{Received February 18, 2013, in f\/inal form August 15, 2013; Published online August 23, 2013}

\Abstract{We consider the space of dif\/ferential operators $\mathcal{D}_{\lambda\mu}$ acting between
$\lambda$- and $\mu$-den\-sities def\/ined on $S^{1|2}$ endowed with its standard contact structure.
This contact structure allows one to def\/ine a~f\/iltration on $\mathcal{D}_{\lambda\mu}$ which is f\/iner
than the classical one, obtained by writting a~dif\/ferential operator in terms of the partial derivatives
with respect to the dif\/ferent coordinates.
The space $\mathcal{D}_{\lambda\mu}$ and the associated graded space of symbols $\mathcal{S}_{\delta}$
($\delta=\mu-\lambda$) can be considered as $\mathfrak{spo}(2|2)$-modules, where $\mathfrak{spo}(2|2)$ is
the Lie superalgebra of contact projective vector f\/ields on $S^{1|2}$.
We show in this paper that there is a~unique isomorphism of $\mathfrak{spo}(2|2)$-modules between
$\mathcal{S}_{\delta}$ and $\mathcal{D}_{\lambda\mu}$ that preserves the principal symbol (i.e.\
an {$\mathfrak{spo}(2|2)$-equivariant} quantization) for some values of $\delta$ called non-critical values.
Moreover, we give an explicit formula for this isomorphism, extending in this way the results
of~[Mellouli~N., \textit{SIGMA} \textbf{5} (2009), 111, 11~pages] which were established for second-order dif\/ferential operators.
The method used here to build the $\mathfrak{spo}(2|2)$-equivariant quantization is the same as the one
used in~[Mathonet P., Radoux F., \textit{Lett. Math. Phys.} \textbf{98}
  (2011), 311--331] to prove the existence of a~$\mathfrak{pgl}(p+1|q)$-equivariant quantization on
$\mathbb{R}^{p|q}$.}


\Classification{53D10; 17B66; 17B10}

\vspace{-2mm}

\section{Introduction}

The concept of equivariant quantization over $\mathbb{R}^{n}$ was introduced by P.~Lecomte and
V.~Ovsienko in~\cite{LO}.
An equivariant quantization is a~linear bijection between a~space of dif\/ferential operators and its
corresponding space of symbols that commutes with the action of a~Lie subalgebra of vector f\/ields over
$\mathbb{R}^{n}$ and preserves the principal symbol.

In their seminal work~\cite{LO}, P.~Lecomte and V.~Ovsienko considered spaces of dif\/ferential
ope\-ra\-tors acting between densities and the Lie algebra of projective vector f\/ields over
$\mathbb{R}^{n}$, \mbox{$\mathfrak{sl}(n+1)$}.
In this situation, they showed the existence and uniqueness of an equivariant quantization.

The results of~\cite{LO} were generalized in many references: in~\cite{BHMP,BM, DLO,Lecras} the authors
considered other spaces of dif\/ferential operators or other Lie subalgebras of vector f\/ields over the
Euclidean space.
In~\cite{Leconj}, P.~Lecomte globalized the problem of equivariant quantization by def\/ining the problem
of natural invariant quantization on arbitrary manifolds.
Finally in~\cite{Bou00,CapSil,Fox,Hansoul,sarah,MR,MR1,MR2,MR3}, the authors proved the existence of such
quantizations by using dif\/ferent methods in more and more general contexts.

Recently, several papers dealt with the problem of equivariant quantizations in the context of
supergeometry: the thesis~\cite{Mic09} dealt with \emph{conformally equivariant quantizations} over
supercotangent bundles, the papers~\cite{MR11} and~\cite{LMR} exposed and solved respectively the problems
of the $\mathfrak{pgl}(p+1|q)$-equivariant quantization over $\mathbb{R}^{p|q}$ and of the
$\mathfrak{osp}(p+1,q+1|2r)$-equivariant quantization over $\mathbb{R}^{p+q|2r}$, whereas in~\cite{LR}, the
authors def\/ine the problem of the natural and projectively invariant quantization on arbitrary
supermanifolds and show the existence of such a~map.

In~\cite{GarMelOvs07,Mel09} the problem of equivariant quantizations over the supercircles $S^{1|1}$ and
$S^{1|2}$ endowed with canonical contact structures was considered.
These quantizations are equivariant with respect to Lie superalgebras of contact projective vector f\/ields.
These Lie superalgebras are the intersections of the Lie superalgebras of contact vector f\/ields on
$S^{1|1}$ and $S^{1|2}$ and of the projective Lie superalgebras $\mathfrak{pgl}(2|1)$ and
$\mathfrak{pgl}(2|2)$.
In these works, the spaces of dif\/ferential operators are endowed with f\/iltrations which are def\/ined
thanks to the contact structures and which are f\/iner than the classical ones.
The spaces of symbols are then the graded spaces corresponding to these f\/iner f\/iltrations.
In~\cite{GarMelOvs07}, the authors show the existence of equivariant quantizations at an arbitrary order
whereas in~\cite{Mel09}, N.~Mellouli proved the existence of equivariant quantizations up to order two.

In this paper, we aim to build an $\mathfrak{spo}(2|2)$-equivariant quantization at an arbitrary order on~$S^{1|2}$, where $\mathfrak{spo}(2|2)$ stands for the Lie superalgebra of contact projective vector
f\/ields on~$S^{1|2}$.
Moreover, we derive an explicit formula for this quantization, extending in this way the results
of~\cite{Mel09}.
The method used here to build the quantization is the same as the one linked to the Casimir operators used
in~\cite{MR11} to build the $\mathfrak{pgl}(p+1|q)$-equivariant quantization on $\mathbb{R}^{p|q}$.

The paper is organized as follows.
In Section~\ref{Section2}, we recall the def\/initions of the objects that occur in the problem of quantization
such as densities, dif\/ferential operators and symbols.
In Section~\ref{Section3}, we expose the tools that we are going to use to build the quantization.
These tools have already been def\/ined in~\cite{MR11}.
The main task performed in Section~\ref{Section3} is the computation and the comparison of the second-order
Casimir operators of $\mathfrak{spo}(2|2)$ acting on the space of dif\/ferential operators and on the space
of symbols on $S^{1|2}$.
Section~\ref{Section4} is devoted to the explicit construction of the $\mathfrak{spo}(2|2)$-equivariant
quantization map, built from the techniques developed in Section~\ref{Section3}.
The method used in Section~\ref{Section4} allows one to f\/ind in Section~\ref{Section5} explicit formulae for the
$\mathfrak{spo}(2|2)$-equivariant quantization at an arbitrary order.

\section{Notation and problem setting}\label{Section2}%\label{tens}

In this section, we recall some tools pertaining to the problem of equivariant quantization such as tensor
(or weighted) densities, dif\/ferential operators, symbols, contact projective vector f\/ields on $S^{1|2}$.
These objects were already exposed in~\cite{Mel09}.

The only point that we will deepen concerns the Lie superalgebra of contact projective vector f\/ields,~$\mathfrak{spo}(2|2)$.
Actually, we will realize this Lie superalgebra as a~Lie subsuperalgebra of~$\mathfrak{gl}(2|2)$, allowing
in this way to def\/ine easily Casimir operators associated with representations of~$\mathfrak{spo}(2|2)$.
This point will be crucial in the sequel.

\subsection[Functions and vector fields on $S^{1|2}$]{Functions and vector f\/ields on $\boldsymbol{S^{1|2}}$}

We def\/ine the supercircle $S^{1|2}$ by describing its graded commutative algebra of functions which we
denote by $C^{\infty}(S^{1|2})$ and which is constituted by the elements
\begin{gather*}
f(x,\theta_{1},\theta_{2})=f_{0}(x)+\theta_{1}f_{1}(x)+\theta_{2}f_{2}(x)+\theta_{1}\theta_{2}f_{12}(x),
\end{gather*}
where $x$ is the coordinate corresponding to one of the two af\/f\/ine coordinates system on $\mathbb{R}
P^{1}$, $\theta_{1}$ and $\theta_{2}$ are odd Grassmann coordinates and where $f_{0},f_{12},f_{1},f_{2}\in
C^{\infty }(S^{1})$ are functions with complex values.
We def\/ine the parity function~$\tilde{\cdot}$ by setting $\tilde{x} =0$ and $\tilde{\theta_{1}}
=\tilde{\theta_{2}}=1$.

A \textit{vector field on} $S^{1|2}$ is a~derivation of the graded commutative algebra
$C^{\infty}(S^{1|2})$.
It can be expressed as
\begin{gather*}
X=f\partial_{x}+g_{1}\partial_{\theta_{1}}+g_{2}\partial_{\theta_{2}},
\end{gather*}
where $f,g_1,g_2\in C^{\infty }(S^{1|2})$, $\partial _{x}=\frac{ \partial }{\partial x}$ and $\partial
_{\theta_{i}}=\frac{\partial }{\partial \theta_{i}}$, for $i=1,2$.
The space of vector f\/ields on~$S^{1|2}$ is a~Lie superalgebra which we shall denote by
$\mathrm{Vect}(S^{1|2})$.

\subsection{The Lie superalgebra of contact vector f\/ields}

The \textit{standard contact} structure on
$S^{1|2}$ is def\/ined by the data of a~linear distribution $\left\langle \overline{D}_{1},\overline{D}
_{2}\right\rangle $ on $S^{1|2}$ generated by the odd vector f\/ields
\begin{gather*}
\overline{D}_{1}=\partial_{\theta_{1}}-\theta_{1}\partial_{x},
\qquad
\overline{D}_{2}=\partial_{\theta_{2}}-\theta_{2}\partial_{x}.
\end{gather*}

A vector f\/ield $X$ on $S^{1|2}$ is called a~\textit{contact vector field} if it preserves the contact
distribution, that is, satisf\/ies the condition:
\begin{gather*}
\left[X,\overline{D}_{1}\right]=\psi_{1_{X}}\overline{D}_{1}+\psi_{2_{X}}\overline{D}_{2},
\qquad
\left[X,\overline{D}_{2}\right]=\phi_{1_{X}}\overline{D}_{1}+\phi_{2_{X}}\overline{D}_{2},
\end{gather*}
where $\psi _{1_{X}},\psi _{2_{X}},\phi _{1_{X}},\phi _{2_{X}}\in C^{\infty }(S^{1|2})$ are functions
depending on $X$.
The space of contact vector f\/ields is a~Lie superalgebra which we shall denote by $\mathcal{K}(2)$.

It is well-known that every contact vector f\/ield can be expressed, for some function $f\in
C^{\infty}(S^{1|2})$, by
\begin{gather*}
X_{f}=f\partial_{x}-\left(-1\right)^{\tilde{f}}\tfrac{1}{2}\left(\overline{D}_{1}\left(f\right)\overline{D}
_{1}+\overline{D}_{2}\left(f\right)\overline{D}_{2}\right).
\end{gather*}

The function $f$ is said to be a~\textit{contact Hamiltonian} of the f\/ield $ X_{f}$.
The space $C^{\infty}(S^{1|2})$ is therefore identif\/ied with the Lie superalgebra $\mathcal{\;K}(2)$ and
is equipped with a~structure of Lie superalgebra thanks to the following contact bracket:
\begin{gather*}
\left\{f,g\right\}=fg^{\prime}-f^{\prime}g-\left(-1\right)^{\tilde{f}}\tfrac{1}{2}\left(\overline{D}_{1}
\left(f\right)\overline{D}_{1}\left(g\right)+\overline{D}_{2}\left(f\right)\overline{D}_{2}
\left(g\right)\right),
\end{gather*}
where $f^{\prime }=\partial _{x}(f)$.

\subsection[The Lie superalgebra $\mathfrak{spo}(2|2)$]{The Lie superalgebra $\boldsymbol{\mathfrak{spo}(2|2)}$}\label{spo}

The Lie superalgebra $\mathfrak{spo}(2|2)$ is the intersection of the Lie superalgebra $\mathcal{\;K}(2)$
and the Lie superalgebra of projective vector f\/ields $\mathfrak{pgl}(2|2)$ exposed in~\cite{MR11}.
The Lie superalgebra $\mathfrak{spo}(2|2)$ is thus a~$4|4$-dimensional Lie superalgebra spanned by the
contact vector f\/ields associated with the following contact Hamiltonians:
\begin{gather*}
\left\{ 1,x,\theta_1,\theta_2,\theta_1\theta_2,x^2,x\theta_1, x\theta_2 \right\}.
\end{gather*}

The Lie subsuperalgebra $\mathfrak{Aff}\left(2|2\right)$ of $\mathfrak{spo} \left(2|2\right)$ spanned by
the contact vector f\/ields associated with the contact Hamiltonians $\left\{
1,x,\theta_1,\theta_2,\theta_1\theta_2\right\}$ will be called the \textit{affine} Lie superalgebra.

Actually, the Lie superalgebra $\mathfrak{spo}(2|2)$ can be realized as the embedding of a~Lie superalgebra
consisted of matrices belonging to $\mathfrak{gl}(2|2)$ into $\mathrm{Vect}(S^{1|2})$.

We shall also denote this matrix realization by $\mathfrak{spo}(2|2)$ (remark that this matrix realization
has not to be confused with the special Poisson Lie superalgebra).
The matrix realization of $\mathfrak{spo}(2|2)$ is e.g.\
exposed in~\cite[p.~419]{LPS}.
In this reference, this Lie superalgebra is denoted by $\mathfrak{osp}^{sk}(2|2)$; this is the Lie
subsuperalgebra of $\mathfrak{gl}(2|2)$ made of the matrices $A$ that preserve a~particular
superskewsymmetric even bilinear form $\omega$ def\/ined on $\mathbb{R}^{2|2}$, in the sense that
\begin{gather*}
\omega(AU,V)+(-1)^{\tilde{A}\tilde{U}}\omega(U,AV)=0
\qquad
\mbox{for all}
\quad
U,V\in\mathbb{R}^{2|2}.
\end{gather*}
This particular form $\omega$ is def\/ined on $\mathbb{R}^{2|2}$ by $\omega(U,V)=V^tGU$, where
\begin{gather*}
G=\left(
\begin{matrix}
J&0
\\
0&S
\end{matrix}
\right),
\qquad
J=\left(
\begin{matrix}
0&-1
\\
1&0
\end{matrix}
\right),
\qquad
S=\left(
\begin{matrix}
1&0
\\
0&1
\end{matrix}
\right).
\end{gather*}

The Lie superalgebra $\mathfrak{spo}(2|2)$ is then constituted by the matrices $A$ in $\mathfrak{gl}(2|2)$
such that
\begin{gather*}
A^{st}G+GA=0,
\end{gather*}
where the supertranspose of the matrix $A$, $A^{st}$, is def\/ined by
\begin{gather*}
\left(
\begin{matrix}
A_{1}^{t}&-A_{3}^{t}
\\
A_{2}^{t}&A_{4}^{t}
\end{matrix}
\right)
\end{gather*}
if the matrix $A$ is equal to
\begin{gather*}
\left(
\begin{matrix}
A_{1}& A_{2}
\\
A_{3}&A_{4}
\end{matrix}
\right).
\end{gather*}

The blocks constituting the matrix $A$ have thus to verify the following properties:
\begin{itemize}\itemsep=0pt
\item
{}$A_{1}^{t}J+JA_{1}=0$, i.e.\
{}$A_{1}\in\mathfrak{sp}(2)$;
\item
{}$A_{4}^{t}+A_{4}=0$, i.e.\
{}$A_{4}\in\mathfrak{o}(2)$;
\item
{}$A_{3}=-A_{2}^{t}J$.
\end{itemize}

The Lie superalgebra $\mathfrak{spo}(2|2)$ is then $4|4$-dimensional and one of its bases is constituted by
the following matrices:
\begin{gather*}
\left(
\begin{matrix}
1&0&0&0
\\
0&-1&0&0
\\
0&0&0&0
\\
0&0&0&0
\\
\end{matrix}
\right),\qquad
\left(
\begin{matrix}
0&1&0&0
\\
0&0&0&0
\\
0&0&0&0
\\
0&0&0&0
\\
\end{matrix}
\right),\qquad
\left(
\begin{matrix}
0&0&0&0
\\
1&0&0&0
\\
0&0&0&0
\\
0&0&0&0
\\
\end{matrix}
\right),
\qquad
\left(
\begin{matrix}
0&0&0&0
\\
0&0&0&0
\\
0&0&0&1
\\
0&0&-1&0
\\
\end{matrix}
\right),
\\
\left(
\begin{matrix}
0&0&1&0
\\
0&0&0&0
\\
0&1&0&0
\\
0&0&0&0
\\
\end{matrix}
\right),\qquad
\left(
\begin{matrix}
0&0&0&1
\\
0&0&0&0
\\
0&0&0&0
\\
0&1&0&0
\\
\end{matrix}
\right),\qquad
\left(
\begin{matrix}
0&0&0&0
\\
0&0&1&0
\\
-1&0&0&0
\\
0&0&0&0
\\
\end{matrix}
\right), \qquad
\left(
\begin{matrix}
0&0&0&0
\\
0&0&0&1
\\
0&0&0&0
\\
-1&0&0&0
\\
\end{matrix}
\right).
\end{gather*}

The upper-left blocks of the f\/irst three matrices provide a~basis of $\mathfrak{sp}(2)$, whereas the
lower-right block of the fourth provides a~basis of $\mathfrak{o}(2)$.

With the above considerations, it is obvious that $\mathfrak{spo}(2|2)$ can be embedded into the projective
Lie superalgebra $\mathfrak{pgl}(2|2)$ via the map $\iota$ def\/ined in the following way:
\begin{gather*}
\iota: \ \mathfrak{spo}(2|2)\to\mathfrak{pgl}(2|2): \ A\mapsto [A].
\end{gather*}

Now, $\mathfrak{pgl}(2|2)$ can be embedded into the Lie superalgebra of vector f\/ields on $S^{1|2}$ thanks
to the projective embedding def\/ined in~\cite{MR11} in the following way:
\begin{gather*}
\left[\left(
\begin{matrix}
0&\xi
\\
v&B
\end{matrix}
\right)\right]\mapsto-\sum_{i=1}^{3}v^i\partial_{y^i}-\sum_{i,j=1}^{3}(-1)^{\tilde{j}(\tilde{i}+\tilde{j})}
B_j^i y^j\partial_{y^i}+\sum_{j=1}^{3}(-1)^{\tilde{j}}\xi_j y^j y^i\partial_{y^i},
\end{gather*}
where $v\in \mathbb{R}^{1|2}$, $\xi\in \mathbb{R}^{1|2*}$, $B\in\mathfrak{gl}(1|2)$ and where the
coordinates $y^{1}$, $y^{2}$, $y^{3}$ correspond respectively to the coordinates $x$, $\theta_{1}$, $\theta_{2}$.

Composing $\iota$ with the projective embedding, we can embed $\mathfrak{spo}(2|2)$ into the Lie
superalgebra of vector f\/ields on $S^{1|2}$.
If we compute this embedding on the generators of $\mathfrak{spo}(2|2)$ written above, we obtain
respectively $2X_{x}$, $X_{x^2}$, $-X_{1}$, $2X_{\theta_{1}\theta_{2}}$, $-2X_{x\theta_{1}}$,
$-2X_{x\theta_{2}}$, $2X_{\theta_1}$ and $2X_{\theta_2}$.

\subsection{Modules of weighted densities}

For any contact vector f\/ield $X_{f}$, we def\/ine a~family of dif\/ferential operators of order one on
$C^{\infty}(S^{1|2})$, denoted by $L_{X_{f}}^{\lambda}$, in the following way:
\begin{gather*}
L_{X_{f}}^{\lambda}=X_{f}+\lambda f^{\prime},
\end{gather*}
where the parameter $\lambda$ is an arbitrary (complex) number and where $f^{\prime}$ denotes the left
multiplication by $f^{\prime}$.
The map $X_{f}\mapsto L_{X_{f}}^{\lambda }$ is a~homomorphism of Lie superalgebras.
We thus obtain a~family of $\mathcal{K}\left( 2\right) $-modules on $C^{\infty}(S^{1|2}) $ which we shall
denote by $\mathcal{F}_{\lambda}$ and which we shall call spaces of \textit{weighted densities} of weight
$\lambda$.

\subsection{Dif\/ferential operators and symbols}%\label{XXX}

In the sequel, we will call ``natural number'' a~non-negative integer.

For every (half)-natural number $k$, we denote by $\mathcal{D}_{\lambda\mu }^{k}$ the space of
dif\/ferential operators acting between $\lambda$- and $\mu$-densities that are of the form
\begin{gather*}
\sum_{l+\frac{m}{2}+\frac{n}{2}\leq k}a_{l,m,n}\left(\partial_{x}\right)^{l}\overline{D}_{1}^{m}\overline{D}
_{2}^{n},
\end{gather*}
where $a_{l,m,n}\in C^{\infty}(S^{1|2})$ for all $l$, $m$, $n$.
Furthermore, since $\partial _{x}=-\overline{D}_{1}^{2}=-\overline{D} _{2}^{2} $, we can assume $m,n\leq 1$.
The space $\mathcal{D}_{\lambda\mu }^{k}$ will be called the space of dif\/ferential operators of order $k$.

We def\/ine then $\mathcal{D}_{\lambda\mu}$, the space of dif\/ferential operators acting between
$\lambda$- and $\mu$-densities as the union of the spaces $\mathcal{D}_{\lambda\mu }^{k}$:
\begin{gather*}
\mathcal{D}_{\lambda\mu}=\bigcup_{k\in\frac {\mathbb{N}}{2}}\mathcal{D}_{\lambda\mu }^{k},
\end{gather*}
where $\frac {\mathbb{N}}{2}$ denotes the set consisting of natural and half-natural numbers.
Since $\mathcal{D}_{\lambda \mu}^{k}\subset\mathcal{D}_{\lambda\mu}^{k+\frac {1}{2}}$, the space
$\mathcal{D}_{\lambda\mu}$ is a~f\/iltered space.
This space has a~structure of a~$\mathcal{K}\left( 2\right)$-module def\/ined in the following way: if
$D\in\mathcal{D}_{\lambda\mu}$, then the Lie derivative of $D$ in the direction of $X_{f}$, denoted by~$\mathcal{L}_{X_f}D$, is given by the dif\/ferential operator
\begin{gather*}
L_{X_f}^{\mu}\circ D-(-1)^{\tilde{f}\tilde{D}}D\circ L_{X_f}^{\lambda}.
\end{gather*}
Moreover, it turns out that the action of $\mathcal{K}(2)$ preserves the order of~$D$.

The graded space associated with the f\/iltered space $\mathcal{D}_{\lambda\mu}$ is called the space of
symbols and is denoted by $\mathcal{S}_{\delta}$, where $\delta=\mu-\lambda$:
\begin{gather*}
\mathcal{S}_{\delta}=\bigoplus_{k\in\frac {\mathbb{N}}{2}} \mathcal{S}_{\delta}^{k},
\end{gather*}
where $\mathcal{S}_{\delta}^{k}=\mathcal{D}_{\lambda \mu }^{k}/\mathcal{D}_{\lambda\mu }^{k-\frac {1}{2}}$
for every (half)-natural number $k$.

The principal symbol map, denoted by $\sigma$, is the map def\/ined on $\mathcal{D}_{\lambda\mu}$ whose the
restriction to $\mathcal{D}_{\lambda\mu}^{k}$, denoted by $\sigma_{k}$, is simply def\/ined by:
\begin{gather*}
\sigma_{k}: \ \mathcal{D}_{\lambda\mu}^{k}\to\mathcal{S}_{\delta}^{k}: \ D\mapsto [D].
\end{gather*}

If $k$ is a~natural number, the space of symbols of degree $k$, $\mathcal{S}_{\delta}^{k}$, is isomorphic
as vector space to $\mathcal{F}_{\delta-k}\oplus\mathcal{F}_{\delta-k}$ through the following
identif\/ication:
\begin{gather*}
\big[F_{1}\partial_{x}^{k}+F_{2}\partial_{x}^{k-1}\bar{D}_{1}\bar{D}_{2}\big]\longleftrightarrow(F_{1},F_{2}).
\end{gather*}

Thanks to the fact that the action of $\mathcal{K}(2)$ on $\mathcal{D}_{\lambda\mu}$ preserves the
f\/iltration of this space, the $\mathcal{K}(2)$-module structure on $\mathcal{D}_{\lambda\mu}^{k}$ induces
a~$\mathcal{K}(2)$-module structure on $\mathcal{S}_{\delta}^{k}$.
If we denote by $L_{X_f}$ the Lie derivative of a~symbol in the direction of $X_{f}$, we have
\begin{gather*}
L_{X_f}[D]:=[\mathcal{L}_{X_f}D].
\end{gather*}

If a~symbol of degree $k$ is represented by a~pair of densities $(F_1,F_2)$, it is easy to see that
$L_{X_f}(F_1,F_2)$ corresponds to the pair
\begin{gather*}
\big(L_{X_f}^{\delta-k}F_{1},L_{X_f}^{\delta-k}F_{2}\big).
\end{gather*}

If $k$ is half of a~natural number, the space of symbols of degree $k$, $\mathcal{S}_{\delta}^{k}$, is also
isomorphic as vector space to $\mathcal{F}_{\delta-k}\oplus\mathcal{F}_{\delta-k}$ through the following
identif\/ication:
\begin{gather*}
\Big[F_{1}\partial_{x}^{k-\frac {1}{2}}\bar{D}_{1}+F_{2}\partial_{x}^{k-\frac{1}{2}}\bar{D}_{2}\Big]
\longleftrightarrow(F_{1},F_{2}).
\end{gather*}

As in the case where $k$ is a~natural number, the $\mathcal{K}(2)$-module structure on
$\mathcal{D}_{\lambda\mu}^{k}$ induces a~$\mathcal{K}(2)$-module structure on $\mathcal{S}_{\delta}^{k}$.
If a~symbol of degree $k$ is represented by a~pair of densities $(F_1,F_2)$, it is easy to see that
$L_{X_f}(F_1,F_2)$ corresponds to the pair
\begin{gather*}
\left(L_{X_f}^{\delta-k}F_{1}-\frac {1}{2}\bar{D}_{1}\bar{D}_{2}(f)F_{2},L_{X_f}^{\delta-k}F_{2}
+\frac{1}{2}\bar{D}_{1}\bar{D}_{2}(f)F_{1}\right).
\end{gather*}

We are now in position to def\/ine the notion of $\mathfrak{spo}(2|2)$-equivariant quantization:
\begin{definition}\label{quanti}
An $\mathfrak{spo}(2|2)$-equivariant quantization is a~linear bijection $Q$ between $\mathcal{S}_{\delta}$
and $\mathcal{D}_{\lambda\mu}$ that preserves the principal symbol in the sense that $\sigma\circ Q={\rm
Id}$ and that commutes with the action of $\mathfrak{spo}(2|2)$, in the sense that
\begin{gather*}
\mathcal{L}_{X_f}\circ Q=Q\circ L_{X_{f}}
\end{gather*}
for all $X_{f}\in\mathfrak{spo}(2|2)$.
\end{definition}

\section{Tools used to build the quantization}\label{Section3}%\label{cons}

In order to tackle the problem of $\mathfrak{spo}(2|2)$-equivariant quantization, we will need to adapt the
tools used in~\cite{MR11} for the case where $\mathfrak{g}= \mathfrak{pgl}(p+1|q)$.
The main ingredients are the af\/f\/ine quantization map and the dif\/ference between the representations
$(\mathcal{S}_{\delta},\mathcal{L})$ and $(\mathcal{S}_{\delta},L)$ of $\mathfrak{spo}(2|2)$ measured by
the map $\gamma$ (see Section~\ref{gamma}).
In the sequel, $\frac {1}{2}+\mathbb{N}$ will denote the set of half-natural numbers.

\subsection{The af\/f\/ine quantization map}

The \emph{affine quantization map} $Q_{\mathrm{Af\/f}}$ is the linear bijection between
$\mathcal{S}_{\delta}$ and $\mathcal{D}_{\lambda\mu}$ def\/ined in the following way:
\begin{gather*}
{\left.Q_{\rm Af\/f}\right|}_{\mathcal{S}_\delta^k}(F_1,F_2)=
\begin{cases}
F_{1}\partial_{x}^k+F_{2}\partial_{x}^{k-1}\bar{D}_1\bar{D}_2,&
\text{if}
\
k\in\mathbb{N},
\\
F_{1}\partial_{x}^{k-\frac{1}{2}}\bar{D}_1+F_{2}\partial_{x}^{k-\frac{1}{2}}\bar{D}_2,&
\text{if}
\
k\in\frac{1}{2}+\mathbb{N}.
\end{cases}
\end{gather*}

Thanks to the map $Q_{\mathrm{Af\/f}}$, we can carry the ${\rm Vect}(S^{1|2})$-module structure of
$\mathcal{D}_{\lambda,\mu}$ to $\mathcal{S}_{\delta}$ by def\/ining a~Lie derivative $\mathcal{L}$ on
$\mathcal{S}_{\delta}$ in the following way:{\samepage
\begin{gather*}
\mathcal{L}_X = Q_{\mathrm{Af\/f}}^{-1}\circ\mathcal{L}_X\circ Q_{\mathrm{Af\/f}}
\end{gather*}
for all $X$ in ${\rm Vect}(S^{1|2})$.}

The existence of an $\mathfrak{spo}(2|2)$-equivariant quantization is then equivalent to the existence, for
all (half)-natural number $k$, of an $\mathfrak{spo}(2|2)$-equivariant map
\begin{gather*}
Q:(\mathcal{S}_{\delta}^{k},L)\to(\mathcal{S}_{\delta},\mathcal{L})
\end{gather*}
such that the homogeneous part of highest degree of $Q(S)$ is equal to $S$ for all
$S\in\mathcal{S}_{\delta}^{k}$.
Indeed, $Q$ is a~map of this type if and only if $Q_{\mathrm{Af\/f}}\circ Q$ is an
$\mathfrak{spo}(2|2)$-equivariant quantization in the sense of Def\/inition~\ref{quanti}.

\subsection[The map $\gamma$]{The map $\boldsymbol{\gamma}$}\label{gamma}

The dif\/ference between the representations $(\mathcal{S}_{\delta},\mathcal{L})$ and
$(\mathcal{S}_{\delta},L)$ of $\mathfrak{spo}(2|2)$ is measured by the map
\begin{gather*}
\gamma: \ \mathfrak{spo}(2|2) \to \mathfrak{gl}(\mathcal{S}_{\delta},\mathcal{S}_{\delta}): \ X_f \mapsto
\gamma(X_f)=\mathcal{L}_{X_f}-L_{X_f}.
\end{gather*}

In order to compute this map, the following lemma, which gives some commutators, will be useful.
\begin{lemma}\label{formules}
If $k\in\mathbb{N}$ and $X_{f}\in\mathfrak{spo}(2|2)$, then, if we consider that the operators $X_{f}$,
$\partial_{x}$, $\bar{D}_{1}$ and~$\bar{D}_{2}$ act on $\mathcal{F}_{\lambda}$ $($where
$\lambda\in\mathbb{R})$, we obtain the following relations:
\begin{gather*}
[X_{f},\partial_{x}^{k}]=-kf'\partial_{x}^{k}+k\frac
{(-1)^{\tilde{f}}}{2}(\bar{D}_{1}(f')\partial_{x}^{k-1}\bar{D}_{1}
+\bar{D}_{2}(f')\partial_{x}^{k-1}\bar{D}_{2})-\frac {k(k-1)}{2}f''\partial_{x}^{k-1},
\\
[X_{f},\bar{D}_{1}]=-\frac {1}{2}f'\bar{D}_{1}+\frac {1}{2}\bar{D}_{1}\bar{D}_{2}(f)\bar{D}_{2},
\qquad
[X_{f},\bar{D}_{2}]=-\frac {1}{2}f'\bar{D}_{2}-\frac {1}{2}\bar{D}_{1}\bar{D}_{2}(f)\bar{D}_{1}
\end{gather*}
and
\begin{gather*}
[\partial_{x}^{k},f]=kf'\partial_{x}^{k-1}+\frac{k(k-1)}{2}f''\partial_{x}^{k-2},
\end{gather*}
if $f$ denotes the left multiplication by $f$.
\end{lemma}
\begin{proof}
The f\/irst formula can be proved simply by induction.
If $k=1$, the relation is obvious.
If the formula is true for $k$ then
\begin{gather*}
[X_{f},\partial_{x}^{k+1}]=[X_{f},\partial_{x}\partial_{x}^{k}]=\left(-f'\partial_{x}+\frac{(-1)^{\tilde{f}}
}{2}(\bar{D}_{1}(f')\bar{D}_{1}+\bar{D}_{2}(f')\bar{D}_{2})\right)\partial_{x}^{k}
\\
\qquad{}
+\partial_{x}\left(-kf'\partial_{x}^{k}+k\frac{(-1)^{\tilde{f}}}{2}(\bar{D}_{1}(f')\partial_{x}^{k-1}\bar{D}
_{1}+\bar{D}_{2}(f')\partial_{x}^{k-1}\bar{D}_{2})-\frac{k(k-1)}{2}f''\partial_{x}^{k-1}\right).
\end{gather*}
It is then easy to see that the formula is true also for $k+1$.

The formulae for $[X_{f},\bar{D}_{1}]$ and $[X_{f},\bar{D}_{2}]$ can be proved simply by a~straightforward
computation while the last formula can be proved e.g.\
by an easy induction.
\end{proof}

We are now in position to compute the map $\gamma$.
We shall say that a~symbol $(F_1,F_2)$ is even (resp.\
odd) if $\tilde{F_{1}}=\tilde{F_{2}}=0$ (resp.\
{}$\tilde{F_{1}}=\tilde{F_{2}}=1$).
If $(F_{1},F_{2})$ is a~homogeneous symbol, we shall denote by $\tilde{F}$ the number
$\tilde{F_{1}}=\tilde{F_{2}}$.
\begin{proposition}\label{gamma0*}
The map $\gamma$ vanishes on the affine Lie subsuperalgebra $\mathfrak{Aff}\left(2|2\right)$.

Also, for any $f \!\in\! \{x\theta_1,x\theta_2,x^2\}$, $\gamma(f)$ maps $\mathcal{S}^k_\delta$ to
$\mathcal{S}^{k-\frac{1}{2}}_\delta \!\oplus \mathcal{S}^{k-1}_\delta$.
More precisely, if $(F_1,F_2) \!\in\!\mathcal{S}^k_\delta$ is a~homogeneous symbol, we have
\begin{gather*}
(Q_{\rm Af\/f}\circ\gamma(X_f))(F_1,F_2)=(-1)^{\tilde{f}}\left((-1)^{\tilde{F}}\frac{k}{2}\bar{D}
_1(f')F_1+(-1)^{\tilde{F}}\left(\frac{k}{2}+\lambda\right)\bar{D}_2(f')F_2\right)\partial_{x}^{k-1}\bar{D}_1
\\
\qquad
{}+(-1)^{\tilde{f}}\left((-1)^{\tilde{F}}\frac{k}{2}\bar{D}_2(f')F_1-(-1)^{\tilde{F}}\left(\frac{k}{2}
+\lambda\right)\bar{D}_1(f')F_2\right)\partial_{x}^{k-1}\bar{D}_2
\\
\qquad
{}-k\left(\frac{k-1}{2}+\lambda\right)f''F_1\partial_{x}^{k-1}
-(k-1)\left(\frac{k}{2}+\lambda\right)f''F_2\partial_{x}^{k-2}\bar{D}_1\bar{D_2},
\end{gather*}
for $k \in \mathbb{N}$ and
\begin{gather*}
(Q_{\rm Af\/f}\circ\gamma(X_f))(F_1,F_2)=-(-1)^{\tilde{F}+\tilde{f}}\left(\frac{k-\frac{1}{2}}{2}
+\lambda\right)\left(\bar{D}_1(f')F_1+\bar{D}_2(f')F_2\right)\partial_{x}^{k-\frac{1}{2}}
\\
{}+(-1)^{\tilde{f}}\frac{k-\frac{1}{2}}{2}\left((-1)^{\tilde{F}}\bar{D}_1(f')F_2-(-1)^{\tilde{F}}\bar{D}
_2(f')F_1\right)\partial_{x}^{k-\frac{3}{2}}\bar{D}_1\bar{D_2}
\\
{}+\left(k-\frac{1}{2}\right)\!\left(\left(\frac{-k+\frac{1}{2}}{2}-\lambda\right)f''F_1\right)\partial_{x}
^{k-\frac{3}{2}}\bar{D}_1
+\left(k-\frac{1}{2}\right)\!\left(\left(\frac{-k+\frac{1}{2}}{2}-\lambda\right)f''F_2\right)\partial_{x}
^{k-\frac{3}{2}}\bar{D}_2
\end{gather*}
for $k \in \frac{1}{2}+\mathbb{N}$.
\end{proposition}

\begin{proof}
If $k\in\mathbb{N}$, $Q_{\rm Af\/f}(F_{1},F_{2})$ can be expressed formally in the following matrix form:
\begin{gather*}
Q_{\rm Af\/f}(F_{1},F_{2})=(F_{1},F_{2})\left(
\begin{matrix}
\partial_{x}^{k}
\\
\partial_{x}^{k-1}\bar{D}_{1}\bar{D}_{2}
\end{matrix}
\right).
\end{gather*}
If we denote by $D$ the column vector $\left(
\begin{matrix}
\partial_{x}^{k}
\\
\partial_{x}^{k-1}\bar{D}_{1}\bar{D}_{2}
\end{matrix}
\right)$, $\mathcal{L}_{X_{f}}(Q_{\rm Af\/f}(F_{1},F_{2}))$ is then equal to
\begin{gather*}
(X_{f}+\mu f')((F_1,F_2)D)-(-1)^{\tilde{f}\tilde{F}}(F_{1},F_{2})D(X_{f}+\lambda f'),
\end{gather*}
which is equal to
\begin{gather*}
(X_{f}\cdot F_{1},X_{f}\cdot F_{2})D+(-1)^{\tilde{f}\tilde{F}}(F_{1},F_{2})[X_{f},D]
+(\mu-\lambda)f'(F_{1},F_{2})D
\\
\qquad
{}-(-1)^{\tilde{f}\tilde{F}}\lambda(F_1,F_2)[D,f'].
\end{gather*}
Thanks to the relations established in Lemma~\ref{formules}, we can see that $[X_{f},D]$ is equal to
\begin{gather*}
\left(
\begin{matrix}
-kf'\partial_{x}^{k}+\frac{k}{2}(-1)^{\tilde{f}}\big(\bar{D}_{1}(f')\partial_{x}^{k-1}\bar{D}_{1}+\bar{D}_{2}
(f')\partial_{x}^{k-1}\bar{D}_{2}\big)
\vspace{1mm}\\
-kf'\partial_{x}^{k-1}\bar{D}_{1}\bar{D}_{2}+\frac{k}{2}(-1)^{\tilde{f}}\big(\bar{D}_{2}(f')\partial_{x}^{k-1}
\bar{D}_{1}-\bar{D}_{1}(f')\partial_{x}^{k-1}\bar{D}_{2}\big)
\end{matrix}
\right)
\\
\qquad
{}+\left(
\begin{matrix}
-\frac{k(k-1)}{2}f''\partial_{x}^{k-1}
\vspace{1mm}\\
-\frac{k(k-1)}{2}f''\partial_{x}^{k-2}\bar{D}_{1}\bar{D}_{2}
\end{matrix}
\right),
\end{gather*}
while $[D,f']$ is equal to
\begin{gather*}
\left(
\begin{matrix}
kf''\partial_{x}^{k-1}
\vspace{1mm}\\
(-1)^{\tilde{f}+1}(\bar{D}_{2}f')\partial_{x}^{k-1}\bar{D}_{1}+(-1)^{\tilde{f}}(\bar{D}_{1}f')\partial_{x}
^{k-1}\bar{D}_{2}+(k-1)f''\partial_{x}^{k-2}\bar{D}_{1}\bar{D}_{2}
\end{matrix}
\right).
\end{gather*}
The terms $(X_{f}\cdot F_{1},X_{f}\cdot F_{2})D$, $(\mu-\lambda) f'(F_{1},F_{2})D$ and
\begin{gather*}
(-1)^{\tilde{f}\tilde{F}}(F_{1},F_{2})\left(
\begin{matrix}
-kf'\partial_{x}^{k}
\vspace{1mm}\\
-kf'\partial_{x}^{k-1}\bar{D}_{1}\bar{D}_{2}
\end{matrix}
\right)
\end{gather*}
are used to reconstruct $Q_{\rm Af\/f}(L_{X_f}(F_{1},F_{2}))$.
It is easy to see that the sum of the other terms is equal to $Q_{\rm Af\/f}(\gamma(X_f)(F_1,F_2))$.

If $k\in\frac {1}{2}+\mathbb{N}$, $Q_{\rm Af\/f}(F_{1},F_{2})$ can be expressed formally in the following
matrix form:
\begin{gather*}
Q_{\rm Af\/f}(F_{1},F_{2})=(F_{1},F_{2})\left(
\begin{matrix}
\partial_{x}^{k-\frac{1}{2}}\bar{D}_{1}
\\
\partial_{x}^{k-\frac{1}{2}}\bar{D}_{2}
\end{matrix}
\right).
\end{gather*}
If we denote by $D$ the column vector $\left(
\begin{matrix}
\partial_{x}^{k-\frac {1}{2}}\bar{D}_{1}
\\
\partial_{x}^{k-\frac {1}{2}}\bar{D}_{2}
\end{matrix}
\right)$, $\mathcal{L}_{X_{f}}(Q_{\rm Af\/f}(F_{1},F_{2}))$ is then equal to
\begin{gather*}
(X_{f}+\mu f')((F_1,F_2)D)-(-1)^{\tilde{f}(\tilde{F}+1)}(F_{1},F_{2})D(X_{f}+\lambda f'),
\end{gather*}
which is equal to
\begin{gather*}
(X_{f}\cdot F_{1},X_{f}\cdot F_{2})D+(-1)^{\tilde{f}\tilde{F}}(F_{1},F_{2})[X_{f},D]+(\mu-\lambda)f'(F_{1}
,F_{2})D
\\
\qquad
{}-(-1)^{\tilde{f}(\tilde{F}+1)}\lambda(F_1,F_2)[D,f'].
\end{gather*}
Thanks to the relations established in the Lemma~\ref{formules}, we can see that $[X_{f},D]$ is equal to
\begin{gather*}
\left(
\begin{matrix}
-kf'\partial_{x}^{k-\frac{1}{2}}\bar{D}_{1}+\frac{1}{2}\bar{D}_{1}\bar{D}_{2}(f)\partial_{x}^{k-\frac{1}
{2}}\bar{D}_{2}
\vspace{1mm}\\
-kf'\partial_{x}^{k-\frac{1}{2}}\bar{D}_{2}-\frac{1}{2}\bar{D}_{1}\bar{D}_{2}(f)\partial_{x}^{k-\frac{1}{2}}
\bar{D}_{1}
\end{matrix}
\right)
\\
\qquad
{}-\left(
\begin{matrix}
\frac{k-\frac{1}{2}}{2}(-1)^{\tilde{f}}\Big(\bar{D}_{1}(f')\partial_{x}^{k-\frac{1}{2}}+\bar{D}_{2}
(f')\partial_{x}^{k-\frac{3}{2}}\bar{D}_{1}\bar{D}_{2}\Big)+\frac{(k-\frac{1}{2})^{2}}{2}f''\partial_{x}
^{k-\frac{3}{2}}\bar{D}_{1}
\vspace{1mm}\\
-\frac{k-\frac{1}{2}}{2}(-1)^{\tilde{f}}\Big(\bar{D}_{1}(f')\partial_{x}^{k-\frac{3}{2}}\bar{D}_{1}\bar{D}_{2}
-\bar{D}_{2}(f')\partial_{x}^{k-\frac{1}{2}}\Big)+\frac{(k-\frac{1}{2})^{2}}{2}f''\partial_{x}^{k-\frac{3}{2}}
\bar{D}_{2}
\end{matrix}
\right),
\end{gather*}
while $[D,f']$ is equal to
\begin{gather*}
\left(
\begin{matrix}
(-1)^{\tilde{f}}\left(k-\frac{1}{2}\right)f''\partial_{x}^{k-\frac{3}{2}}\bar{D}_{1}+\bar{D}_{1}
(f')\partial_{x}^{k-\frac{1}{2}}
\vspace{1mm}\\
(-1)^{\tilde{f}}\left(k-\frac{1}{2}\right)f''\partial_{x}^{k-\frac{3}{2}}\bar{D}_{2}+\bar{D}_{2}
(f')\partial_{x}^{k-\frac{1}{2}}
\end{matrix}
\right).
\end{gather*}
The terms $(X_{f}\cdot F_{1},X_{f}\cdot F_{2})D$, $(\mu-\lambda) f'(F_{1},F_{2})D$ and
\begin{gather*}
(-1)^{\tilde{f}\tilde{F}}(F_{1},F_{2})\left(
\begin{matrix}
-kf'\partial_{x}^{k-\frac{1}{2}}\bar{D}_{1}+\frac{1}{2}\bar{D}_{1}\bar{D}_{2}(f)\partial_{x}^{k-\frac{1}
{2}}\bar{D}_{2}
\vspace{1mm}\\
-kf'\partial_{x}^{k-\frac{1}{2}}\bar{D}_{2}-\frac{1}{2}\bar{D}_{1}\bar{D}_{2}(f)\partial_{x}^{k-\frac{1}{2}}
\bar{D}_{1}
\end{matrix}
\right)
\end{gather*}
are used to reconstruct $Q_{\rm Af\/f}(L_{X_f}(F_{1},F_{2}))$.
It is easy to see that the sum of the other terms is equal to $Q_{\rm Af\/f}(\gamma(X_f)(F_1,F_2))$.
\end{proof}

If $k\in\mathbb{N}$, then the component of $\gamma(X_{f})|_{\mathcal{S}_{\delta}^{k}}$ with respect to
${\mathcal{S}_{\delta}^{k-\frac {1}{2}}}$ can be put in the following matrix form:
\begin{gather*}
(-1)^{\tilde{f}+\tilde{F}}\left(
\begin{matrix}
\frac{k}{2}\bar{D}_{1}(f')&\big(\frac{k}{2}+\lambda\big)\bar{D}_{2}(f')
\vspace{1mm}\\
\frac{k}{2}\bar{D}_{2}(f')&-\big(\frac{k}{2}+\lambda\big)\bar{D}_{1}(f')
\end{matrix}
\right)
\end{gather*}
while the component with respect to ${\mathcal{S}_{\delta}^{k-1}}$ can be written in the following way:
\begin{gather*}
-f''\left(
\begin{matrix}
k\left(\frac{k-1}{2}+\lambda\right)&0
\\
0&(k-1)\left(\frac{k}{2}+\lambda\right)
\end{matrix}
\right).
\end{gather*}

If $k\in\frac {1}{2}+\mathbb{N}$, then the component of $\gamma(X_{f})|_{\mathcal{S}_{\delta}^{k}}$ with
respect to ${\mathcal{S}_{\delta}^{k-\frac {1}{2}}}$ can be put in the following matrix form:
\begin{gather*}
-(-1)^{\tilde{f}+\tilde{F}}\left(
\begin{matrix}
\left(\frac{k-\frac{1}{2}}{2}+\lambda\right)\bar{D}_{1}(f')&\left(\frac{k-\frac{1}{2}}{2}+\lambda\right)\bar{D}_{2}
(f')
\vspace{1mm}\\
\frac{k-\frac{1}{2}}{2}\bar{D}_{2}(f')&-\left(\frac{k-\frac{1}{2}}{2}\right)\bar{D}_{1}(f')
\end{matrix}
\right)
\end{gather*}
while the component with respect to ${\mathcal{S}_{\delta}^{k-1}}$ can be written in the following way:
\begin{gather*}
\left(k-\frac{1}{2}\right)\left(\frac{-k+\frac{1}{2}}{2}-\lambda\right)f''{\rm Id}.
\end{gather*}

\subsection{Casimir operators}

The method that we are going to use here to build the $\mathfrak{spo}(2|2)$-equivariant quantization is
linked to the Casimir operators, as in~\cite{MR11}.
Actually, the method is based on the comparison of the Casimir operators $C$ and $\mathcal{C}$ of
$\mathfrak{spo}(2|2)$ associated with the representations $L$ and $\mathcal{L}$ on $\mathcal{S}_\delta$.
In this section, we are going to show that the Casimir operator $C$ is diagonalizable and that there is
a~simple relation between $C$ and $\mathcal{C}$.

First of all, let us recall the def\/inition of the Casimir operator associated with a~representation of
a~Lie superalgebra (see e.g.~\cite{Ber87,KacAdvances,Mus97,Pin90,Ser99}).
\begin{definition}\label{casimir operator}
We consider a~Lie superalgebra $\mathfrak{l}$ endowed with an even non-degenerate supersymmetric bilinear
form $K$ and a~representation $(V,\beta)$ of $\mathfrak{l}$.
The Casimir operator $C_\beta$ of $\mathfrak{l}$ associated with $(V,\beta)$ is def\/ined by
\begin{gather*}
C_\beta=\sum_{i}\beta\big(u^*_{i}\big)\beta(u_i),
\end{gather*}
where $u_i$ and $u^*_i$ are $K$-dual bases of $\mathfrak{l}$, in the sense that $K(u_i,u_j^*) =
\delta_{i,j}$ for all $i$, $j$.
\end{definition}

In the sequel, the bilinear form that we will use to def\/ine the Casimir operators of
$\mathfrak{spo}(2|2)$ will be the form $K$ def\/ined in this way:
\begin{gather*}%\label{Killing}
K(A,B)=2\,\mathrm{str}(AB)
\qquad
\forall\, A,B\in\mathfrak{spo}(2|2).
\end{gather*}

The following lemma gives the bases that will be used to def\/ine the Casimir operators of
$\mathfrak{spo}(2|2)$.
\begin{lemma}\label{lem:Killing dual bases}
The $K$-dual basis corresponding to the basis
\begin{gather*}
\{X_1,X_{\theta_1},X_{\theta_2},X_x,X_{\theta_1\theta_2},X_{x\theta_1},X_{x\theta_2},X_{x^2}\}
\end{gather*}
of $\mathfrak{spo}(2|2)$ is given by the basis
\begin{gather*}
\left\lbrace
-\frac{1}{2}X_{x^2},-X_{x\theta_1},-X_{x\theta_2},X_{x},
X_{\theta_1\theta_2},X_{\theta_1},X_{\theta_2},-\frac {1}{2}X_{1}\right\rbrace.
\end{gather*}
\end{lemma}
\begin{proof}
The result is easily proved using the correspondence between vector f\/ields of $\mathfrak{spo}(2|2)$ and
matrices established in Section~\ref{spo}.
\end{proof}

We are now in position to compute the Casimir operator associated with the representa\-tion~$(\mathcal{S}^k_\delta,L)$.
Actually, this operator is simply a~multiple of the identity.
\begin{proposition}\label{Cassl}
The Casimir operator associated with the representa\-tion~$(\mathcal{S}^k_\delta,L)$ of $\mathfrak{spo}(2|2)$
is given by $C|_{\mathcal{S}^k_\delta}=\alpha_{k}{\rm Id}$, where
\begin{gather*}
\alpha_{k}=
\begin{cases}
(-k+\delta)^2,
&
\text{if}
\
k\in\mathbb{N},
\\
(-k+\delta)^2-\frac{1}{4},
&
\text{if}
\
k\in\frac{1}{2}+\mathbb{N}.
\end{cases}
\end{gather*}
\end{proposition}
\begin{proof}
First, using Def\/inition~\ref{casimir operator} and Lemma~\ref{lem:Killing dual bases}, the Casimir
operator $C|_{\mathcal{S}^k_\delta}$ is equal to
\begin{gather*} %\begin{gather*}
-\frac{1}{2}L_{X_{x^2}}L_{X_1}-L_{X_{x\theta_1}}L_{X_{\theta_1}}-L_{X_{x\theta_2}}L_{X_{\theta_2}}
\\
\qquad
{}+(L_{X_{x}})^{2}+(L_{X_{\theta_{1}\theta_{2}}})^{2}+L_{X_{\theta_1}}L_{X_{x\theta_1}}
+L_{X_{\theta_2}}L_{X_{x\theta_2}}-\frac {1}{2}L_{X_{1}}L_{X_{x^2}}.
\end{gather*}
If $k\in\mathbb{N}$, thanks to the form of the hamiltonian vector f\/ields, it is easy to see that this
Casimir operator is a~dif\/ferential operator that can be written in the following way:
\begin{gather*}
f_{0}+f_{1}\partial_{x}+f_{2}\bar{D}_{1}+f_{3}\bar{D}_{2}+f_{4}\partial_{x}^{2}+f_{5}\partial_{x}\bar{D}_{1}
+f_{6}\partial_{x}\bar{D}_{2}+f_{7}\bar{D}_{1}\bar{D}_{2},
\end{gather*}
where $f_{i}\in C^{\infty}(S^{1|2})$ for all $i$.
Since the Casimir operator commutes with the vector f\/ield $X_{1}=\partial_{x}$ and since $\partial_{x}$
commutes with $\bar{D}_{1}$ and $\bar{D}_{2}$, the coef\/f\/icients $f_{i}$ can not depend on the
coordinate $x$.

Since the Casimir operator commutes with the action of the vector f\/ield $X_{\theta_{1}}$, which is equal
to $\frac {1}{2}(\partial_{\theta_{1}}+\theta_{1}\partial_{x})$, and the action of the vector f\/ield
$X_{\theta_{2}}$, which is equal to $\frac {1}{2}(\partial_{\theta_{2}}+\theta_{2}\partial_{x})$, and since
these vector f\/ields commute with $\partial_{x}$, $\bar{D}_{1}$ and $\bar{D}_{2}$, the coef\/f\/icients
$f_{i}$ are invariant under the actions of $X_{\theta_{1}}$ and $X_{\theta_{2}}$, thus these
coef\/f\/icients can not depend on the coordinates $\theta_{1}$ and $\theta_{2}$.

To summarize, $C|_{\mathcal{S}^k_\delta}$ has constant coef\/f\/icients.

Because of the expressions of the Lie derivative and of the vector f\/ields $X_{x^2}$, $X_{x\theta_1}$,
$X_{x\theta_2}$ and $X_{\theta_{1}\theta_{2}}$, the terms $-\frac {1}{2}L_{X_{x^2}}L_{X_1}$,
$L_{X_{x\theta_1}}L_{X_{\theta_1}}$, $L_{X_{x\theta_2}}L_{X_{\theta_2}}$ and
$(L_{X_{\theta_{1}\theta_{2}}})^{2}$ of $C|_{\mathcal{S}^k_\delta}$ give then no contribution.

The term $-\frac {1}{2}L_{X_{1}}L_{X_{x^2}}$ gives the same contribution as the term
\begin{gather*}
-\frac{1}{2}L_{[X_{1},X_{x^2}]}=-\frac{1}{2}L_{X_{\lbrace1,x^{2}\rbrace}}=-L_{X_{x}}.
\end{gather*}
The contribution of this last term is equal to $-(\delta-k){\rm Id}$.

The term $(L_{X_{x}})^{2}$ gives as for it a~contribution equal to $(\delta-k)^{2}{\rm Id}$.

Eventually, each of the two terms $L_{X_{\theta_1}}L_{X_{x\theta_1}}$ and
$L_{X_{\theta_2}}L_{X_{x\theta_2}}$ gives a~contribution equal to $\frac {1}{2}(\delta-k){\rm Id}$.

To conclude, if $k\in\mathbb{N}$, $C|_{\mathcal{S}^k_\delta}=(\delta-k)^{2}{\rm Id}.$

If $k\in\frac {1}{2}+\mathbb{N}$, it is easy to see that the component number 1 (resp.\
2) of $C|_{\mathcal{S}^k_\delta}(F_{1},F_{2})$ can be written in the following way:
\begin{gather*}
\big(f_{0}+f_{1}\partial_{x}+f_{2}\bar{D}_{1}+f_{3}\bar{D}_{2}+f_{4}\partial_{x}^{2}+f_{5}\partial_{x}
\bar{D}_{1}+f_{6}\partial_{x}\bar{D}_{2}+f_{7}\bar{D}_{1}\bar{D}_{2}\big)F_{1}
\\
\qquad
{}+\big(g_{0}+g_{1}\partial_{x}+g_{2}\bar{D}_{1}+g_{3}\bar{D}_{2}\big)F_{2}
\\
\big({\rm resp.\ }\;\big(f_{0}+f_{1}\partial_{x}+f_{2}\bar{D}_{1}+f_{3}\bar{D}_{2}+f_{4}\partial_{x}^{2}+f_{5}
\partial_{x}\bar{D}_{1}+f_{6}\partial_{x}\bar{D}_{2}+f_{7}\bar{D}_{1}\bar{D}_{2}\big)F_{2}
\\
\qquad
{}+\big(g_{0}+g_{1}\partial_{x}+g_{2}\bar{D}_{1}+g_{3}\bar{D}_{2}\big)F_{1}\big),
\end{gather*}
where $f_{i},g_{i}\in C^{\infty}(S^{1|2})$ for all $i$.
As above, thanks to the fact that $C|_{\mathcal{S}^k_\delta}$ commutes with the vector f\/ields $X_{1}$,
$X_{\theta_{1}}$ and $X_{\theta_{2}}$, it turns out that the coef\/f\/icients $f_{i}$ and $g_{i}$ are
constant.

It is then easily seen that the only terms that did not occur in the case $k\in\mathbb{N}$ and that give
contributions in the computation of $C|_{\mathcal{S}^k_\delta}$ come from the term
$(L_{X_{\theta_{1}\theta_{2}}})^{2}$ and give a~contribution equal to $-\frac {1}{4}{\rm Id}$.
\end{proof}

As in~\cite{MR11}, we def\/ine an operator which measures the dif\/ference between the Casimir operators~$C$ and~$\mathcal{C}$.
This map will be called the operator $N$, as in the previous reference.
\begin{definition}%\label{casinil0}
The operator $N$ is def\/ined by
\begin{gather*}%\label{casinil}
N:\ \mathcal{S}_\delta\to\mathcal{S}_\delta: \ S\mapsto\mathcal{C}(S)-C(S).
\end{gather*}
\end{definition}

From the expression of the map $\gamma$, it is easy to deduce an explicit formula for $N$.
\begin{proposition}\label{gamma0}
The operator $N$ maps $\mathcal{S}_{\delta}^{k}$ into $\mathcal{S}_{\delta}^{k-\frac
{1}{2}}\oplus\mathcal{S}_{\delta}^{k-1}$.

If $k\in\mathbb{N}$, the component with respect to $\mathcal{S}_{\delta}^{k-\frac {1}{2}}$ of
$N|_{\mathcal{S}_{\delta}^{k}}$ can be written in a~matrix form in the following way:
\begin{gather*}
-(-1)^{\tilde{F}}\left(
\begin{matrix}
\frac{k}{2}\bar{D}_{1}&(\frac{k}{2}+\lambda)\bar{D}_{2}
\vspace{1mm}\\
\frac{k}{2}\bar{D}_{2}&-(\frac{k}{2}+\lambda)\bar{D}_{1}
\end{matrix}
\right),
\end{gather*}
while the component with respect to $\mathcal{S}_{\delta}^{k-1}$ of this operator is equal to
\begin{gather*}
2\left(
\begin{matrix}
k\left(\frac{k-1}{2}+\lambda\right)\partial_{x}&0
\\
0&(k-1)\left(\frac{k}{2}+\lambda\right)\partial_{x}
\end{matrix}
\right).
\end{gather*}

If $k \in \frac{1}{2}+\mathbb{N}$, the component with respect to $\mathcal{S}_{\delta}^{k-\frac {1}{2}}$ of
$N|_{\mathcal{S}_{\delta}^{k}}$ is equal to
\begin{gather*}
(-1)^{\tilde{F}}\left(
\begin{matrix}
\left(\frac{k-\frac{1}{2}}{2}+\lambda\right)\bar{D}_{1}&\left(\frac{k-\frac{1}{2}}{2}+\lambda\right)\bar{D}_{2}
\\
\frac{k-\frac{1}{2}}{2}\bar{D}_{2}&-\frac{k-\frac{1}{2}}{2}\bar{D}_{1}
\end{matrix}
\right),
\end{gather*}
while the component with respect to $\mathcal{S}_{\delta}^{k-1}$ of this operator can be written in the
following way:
\begin{gather*}
2\left(k-\frac{1}{2}\right)\left(\frac{k-\frac{1}{2}}{2}+\lambda\right)\left(
\begin{matrix}
\partial_{x}&0
\\
0&\partial_{x}
\end{matrix}
\right).
\end{gather*}

\end{proposition}
\begin{proof}
Using the expressions of $\mathcal{C}$ and $C$ and the def\/inition of $\gamma$, we obtain immediately
\begin{gather*}
N=-\frac{1}{2}\gamma(X_{x^2})L_{X_1}-\gamma(X_{x\theta_1})L_{X_{\theta_1}}-\gamma(X_{x\theta_2}
)L_{X_{\theta_2}}+L_{X_{\theta_1}}\gamma(X_{x\theta_1})
\\
\hphantom{N=}{}
+L_{X_{\theta_2}}\gamma(X_{x\theta_2})-\frac{1}{2}L_{X_{1}}\gamma(X_{x^2}).
\end{gather*}
Suppose f\/irst that $k\in\mathbb{N}$.
If we write the Lie derivatives in a~matrix way, we obtain then that the component with respect to
$\mathcal{S}_{\delta}^{k-\frac {1}{2}}$ of $N|_{\mathcal{S}^k_\delta}$ is equal to
\begin{gather*}
2(-1)^{\tilde{F}}\!\left(
\begin{matrix}
\frac{k}{2}\theta_1&\big(\frac{k}{2}+\lambda\big)\theta_2
\vspace{1mm}\\
\frac{k}{2}\theta_2&-\big(\frac{k}{2}+\lambda\big)\theta_1
\end{matrix}
\right)\!\left(
\begin{matrix}
\partial_{x}&0
\\
0&\partial_{x}
\end{matrix}
\right)
 -(-1)^{\tilde{F}}\!\left(
\begin{matrix}
\frac{k}{2}&0
\\
0&-\left(\frac{k}{2}+\lambda\right)
\end{matrix}
\right)\!\left(
\begin{matrix}
\partial_{\theta_1}+\theta_{1}\partial_{x}\!&0
\\
0&\partial_{\theta_1}+\theta_{1}\partial_{x}
\end{matrix}
\right)\!
\\
\qquad
{}-(-1)^{\tilde{F}}\left(
\begin{matrix}
0&\frac{k}{2}+\lambda
\\
\frac{k}{2}&0
\end{matrix}
\right)\left(
\begin{matrix}
\partial_{\theta_2}+\theta_{2}\partial_{x}&0
\\
0&\partial_{\theta_2}+\theta_{2}\partial_{x}
\end{matrix}
\right),
\end{gather*}
i.e.\
to
\begin{gather*}
=-(-1)^{\tilde{F}}\left(
\begin{matrix}
\frac {k}{2}\bar{D}_{1}& \big(\frac{k}{2}+\lambda\big)\bar{D}_{2}
\vspace{1mm}\\
\frac {k}{2}\bar{D}_{2} & -\big(\frac{k}{2}+\lambda\big)\bar{D}_{1}
\end{matrix}
\right).
\end{gather*}
It turns out that the component with respect to $\mathcal{S}_{\delta}^{k-1}$ of $N|_{\mathcal{S}^k_\delta}$
is equal to
\begin{gather*}
2\left(
\begin{matrix}
k\left(\frac{k-1}{2}+\lambda\right)&0
\\
0&(k-1)\left(\frac{k}{2}+\lambda\right)
\end{matrix}
\right)\left(
\begin{matrix}
\partial_{x}&0
\\
0&\partial_{x}
\end{matrix}
\right).
\end{gather*}
Similar computations allow to show that the formulae for the case $k\in\mathbb{N}+\frac{1}{2}$ are also
correct.
\end{proof}

\section{Construction of the quantization}\label{Section4}

We begin this section with the def\/inition of the critical values of $\delta$.
As in~\cite{MR11}, the set of critical values of $\delta$ is a~set outside of which the
$\mathfrak{spo}(2|2)$-equivariant quantization exists and is unique.

\subsection{Critical values}
\begin{definition}
A value of $\delta$ is \emph{critical} if there exist $k,l\in\frac {\mathbb{N}}{2}$ with $l<k$ such that
$\alpha_{k}=\alpha_{l}$.
\end{definition}
\begin{proposition}
The set of critical values of $\delta$ is given by
\begin{gather*}
\left\{\frac{k+l}{2}\colon k,l\in\mathbb{N},l<k\right\}\cup\left\{\frac{k^2-l^2+\frac{1}{4}}{2(k-l)}
\colon k\in\mathbb{N},l\in\mathbb{N}+\frac{1}{2},l<k\right\}
\\
\qquad{}
\cup\left\{\frac{k^2-l^2-\frac{1}{4}}{2(k-l)}\colon k\in\mathbb{N}+\frac{1}{2},l\in\mathbb{N},l<k\right\}.
\end{gather*}
\end{proposition}
\begin{proof}
The result can be easily proved from a~straightforward computation using Proposition~\ref{Cassl}.\looseness=1   %
\end{proof}

\subsection{The construction} The main result is given by the following theorem.
\begin{theorem}\label{flatex}
If $\delta$ is not critical, there exists a~unique $\mathfrak{spo}(2|2)$-equivariant quantization from
$\mathcal{S}_{\delta}$ to $\mathcal{D}_{\lambda,\mu}$.
\end{theorem}
\begin{proof}
The proof is similar to the one of~\cite{MR11}.
First, remark that for every $S\in\mathcal{S}_{\delta}^{k}$, there exists a~unique eigenvector $\hat{S}$ of
$\mathcal{C}$ with eigenvalue $\alpha_{k}$ such that
\begin{gather*}
\begin{cases}\hat{S}=S_k+S_{k-\frac{1}{2}}+S_{k-1}+\cdots+S_0,
\qquad
S_k=S,
\\
S_l\in\mathcal{S}_{\delta}^{l}
\qquad
\mbox{for all}
\quad
l\leqslant k-\frac{1}{2}.
\end{cases}
\end{gather*}

Indeed, the fact that $\hat{S}$ has to be an eigenvector of $\mathcal{C}$ with eigenvalue $\alpha_{k}$
implies that{\samepage
\begin{gather}\label{flatP}
\begin{cases}(C-\alpha_{k}\text{Id})S_{k-\frac{1}{2}}=-{\rm pr}_{k-\frac{1}{2}}(N(S_{k})),
\\
(C-\alpha_{k}\text{Id})S_{k-l}
=-\big({\rm pr}_{k-l}(N(S_{k-l+\frac{1}{2}}))+{\rm pr}_{k-l}(N(S_{k-l+1}))\big) %\mathoperator ?
\\
\end{cases}
\end{gather}
for $l=1,\frac {3}{2},\ldots,k$, where ${\rm pr}_{i}$ stands for the natural projection $\mathcal{S}_\delta
\to \mathcal{S}_\delta^i$.}

As $\delta$ is not critical, the dif\/ferences $\alpha_{k}-\alpha_{l}$ are dif\/ferent from $0$, the
operators $(C-\alpha_{k}{\rm Id})\vert_{\mathcal{S}_{\delta}^{k-l}}$ are thus all invertible and therefore
this system of equations has a~unique solution.

Now, def\/ine the quantization $Q$ by
$
Q\vert_{\mathcal{S}_{\delta}^{k}}(S)=\hat{S}$.
It is clearly a~bijection and it also fulf\/ills
\begin{gather*}
Q\circ L_{X_f}=\mathcal{L}_{X_f}\circ Q
\qquad
\mbox{for all}
\quad
X_f\in \mathfrak{spo}(2|2).
\end{gather*}
Indeed, for all $S\in\mathcal{S}_{\delta}^{k}$, the symbols $Q(L_{X_f}S)$ and $ \mathcal{L}_{X_f}(Q(S))$
share the following properties:
\begin{itemize}\itemsep=0pt
\item
they are eigenvectors of $\mathcal{C}$ of eigenvalue $\alpha_{k}$ because, on the one hand, $\mathcal{C}$
commutes with $\mathcal{L}_{X_f}$ for all $X_f\in\mathfrak{spo}(2|2)$ and, on the other hand, $C$ commutes
with $L_{X_f}$ for all $X_f\in\mathfrak{spo}(2|2)$.
\item
their term of degree $k$ is exactly $L_{X_f}S$.
\end{itemize}
The f\/irst part of the proof shows that they have to coincide.
An $\mathfrak{spo}(2|2)$-equivariant quantization is thus given by $Q_{\rm Af\/f}\circ Q$.

This quantization is the unique $\mathfrak{spo}(2|2)$-equivariant quantization.
Indeed, if $S$ is an eigenvector of $C$ of eigenvalue $\alpha_{k}$, then, if $Q$ is
$\mathfrak{spo}(2|2)$-equivariant, $Q(S)$ has to be an eigenvector of $\mathcal{C}$ of the same eigenvalue
because $\mathcal{C}$ and $C$ are built from the Lie derivatives $\mathcal{L}$ and $L$.
The uniqueness of the solution of the system~\eqref{flatP} implies then the uniqueness of the
$\mathfrak{spo}(2|2)$-equivariant quantization.
\end{proof}

\section[Explicit formulae for the $\mathfrak{spo}(2|2)$-equivariant quantization]
{Explicit formulae for the $\boldsymbol{\mathfrak{spo}(2|2)}$-equivariant quantization}\label{Section5}

In this section, we give in the non-critical situations the explicit formula for the
$\mathfrak{spo}(2|2)$-equivariant quantization in the cases where $k\in\mathbb{N}$ and where
$k\in\mathbb{N}+\frac{1}{2}$.

\subsection[Case $k\in\mathbb{N}$]{Case $\boldsymbol{k\in\mathbb{N}}$}

\begin{proposition}\label{expl1}
If $\delta$ is non-critical and if $l$ is a~natural number, the symbol $S_{k-l}$ in the proof of
Theorem~{\rm \ref{flatex}} is given by:
\begin{gather}\label{sk-l1}
C_l
\begin{pmatrix}
A_kB_{k-l}\partial_x^l&0
\\
0&A_{k-l}B_k\partial_x^l
\end{pmatrix}
S_{k}+D_l
\begin{pmatrix}
A_kB_{k-l}\partial_x^l&-B_kB_{k-l}\bar{D}_1\bar{D}_2\partial_x^{l-1}
\\
A_kA_{k-l}\bar{D}_1\bar{D}_2\partial_x^{l-1}&A_{k-l}B_k\partial_x^l
\end{pmatrix}
S_k,
\end{gather}
where the coefficients $A_{k}$, $B_{k}$, $C_{l}$ and $D_{l}$ are given in the following way:
\begin{gather*}
A_k=-k,
\qquad
B_k=-(k+2\lambda),
\\
C_l=\frac{\prod\limits_{i=1}^{l-1}{A_{k-i}B_{k-i}}}{\prod\limits_{i=1}^l{(\alpha_{k}-\alpha_{k-i})}},
\qquad
D_l=\frac{-l\prod\limits_{i=1}^{l-1}{A_{k-i}B_{k-i}}}{2(\alpha_{k}-\alpha_{k-\frac{1}{2}})\prod\limits_{i=1}^l{(\alpha_{k}
-\alpha_{k-i})}}.
\end{gather*}

If $l\in\mathbb{N}+\frac {1}{2}$, we have:
\begin{gather}\label{sk-l2}
S_{k-l}=(-1)^{\tilde{S_{k}}+1}E_l\left[
\begin{pmatrix}
A_k\bar{D_1}&B_k\bar{D}_2
\\
A_k\bar{D}_2&-B_k\bar{D}_1
\end{pmatrix}
\partial_x^{\lfloor l\rfloor}\right]S_k,
\end{gather}
where $\lfloor l \rfloor$ is the integer part of $l$ and where $E_{l}$ is given in the following way:
\begin{gather*}
E_l=\frac{\prod\limits_{i=1}^{l-\frac{1}{2}}{A_{k-i}B_{k-i}}}
{2\Big(\alpha_{k}-\alpha_{k-\frac{1}{2}}\Big)\prod\limits_{i=1}^{l-\frac{1}{2}}{(\alpha_{k}-\alpha_{k-i})}}.
\end{gather*}

The first coefficients $C_l$, $D_l$ and $E_l$ are given by the following formulae:
\begin{gather*}
C_1= \frac{1}{\alpha_{k}-\alpha_{k-1}},
\qquad
D_1=-\frac{1}{2\Big(\alpha_{k}-\alpha_{k-\frac{1}{2}}\Big)(\alpha_{k}-\alpha_{k-1})},
\qquad
E_{\frac{1}{2}}=\frac{1}{2\Big(\alpha_{k}-\alpha_{k-\frac{1}{2}}\Big)}.
\end{gather*}
\end{proposition}
\begin{proof}
The symbols $S_{k-l}$ are built by induction using the formula~\eqref{flatP}.
When $l=\frac {1}{2}$ and when $l=1$, the formulae~\eqref{sk-l1} and~\eqref{sk-l2} can be easily proved
using the relation~\eqref{flatP}, Propositions~\ref{Cassl} and~\ref{gamma0}.

Suppose that the formulae for the symbols $S_{k-l'}$ are true if $l'\leqslant l$ and prove that the formula
for $S_{k-l-\frac {1}{2}}$ is also true.

Using~\eqref{flatP}, we can f\/irst write
\begin{gather}\label{evp2}
S_{k-l-\frac{1}{2}}=\frac{1}{\alpha_k-\alpha_{k-l-\frac{1}{2}}}\left[N^{k-l}_{{k-l-\frac{1}{2}}}(S_{k-l}
)+N^{k-l+\frac{1}{2}}_{k-l-\frac{1}{2}}\Big(S_{k-(l-\frac{1}{2})}\Big)\right],
\end{gather}
where the component with respect to $\mathcal{S}_{\delta}^{k-l-i}$ of $N|_{\mathcal{S}_{\delta}^{k-l}}$ is
denoted by $N_{k-l-i}^{k-l}$.

Two situations may occur:

If $l$ is a~natural number, then $(k-l)$ is a~natural number and $(k-l+\frac{1}{2})\in\mathbb{N}+\frac{1}{2}$.
In this case, the relation~\eqref{evp2} shows, using Propositions~\ref{Cassl} and~\ref{gamma0}, the
relation~\eqref{sk-l1} and the relation~\eqref{sk-l2}, that $S_{k-l-\frac {1}{2}}$ has the same form as the
right hand side of~\eqref{sk-l2} and that
\begin{gather*}
E_{l+\frac{1}{2}}
=\frac{A_{k-l}B_{k-l}}{2\Big(\alpha_k-\alpha_{k-(l+\frac{1}{2})}\Big)}\left(C_l+2D_l+2E_{l-\frac{1}{2}}\right).
\end{gather*}
Using the induction hypothesis on the coef\/f\/icients $C_{l}$, $D_{l}$ and $E_{l-\frac{1}{2}}$, it is easy
to see that $E_{l+\frac {1}{2}}$ is given by the right formula.

If $l\in\mathbb{N}+\frac {1}{2}$, then $(k-l)\in\mathbb{N}+\frac {1}{2}$ and $\big(k-l+\frac{1}{2}\big)$ is
a~natural number.
Computations similar to the computations above show that $S_{k-l-\frac {1}{2}}$ has the same form as the
right hand side of~\eqref{sk-l1} and that the coef\/f\/icients $C_{l+\frac {1}{2}}$ and $D_{l+\frac
{1}{2}}$ are given by the following formulae:
\begin{gather*}
C_{l+\frac{1}{2}}
=\frac{A_{k-l+\frac{1}{2}}B_{k-l+\frac{1}{2}}}{\Big(\alpha_k-\alpha_{k-(l+\frac{1}{2})}\Big)}C_{l-\frac{1}{2}},
\\
D_{l+\frac{1}{2}}
=\frac{1}{\Big(\alpha_k-\alpha_{k-(l+\frac{1}{2})}\Big)}
\left(-E_l+A_{k-l+\frac{1}{2}}B_{k-l+\frac{1}{2}}D_{l-\frac{1}{2}}\right).
\end{gather*}

Using the induction hypothesis on the coef\/f\/icients $C_{l-\frac {1}{2}}$, $E_{l}$ and
$D_{l-\frac{1}{2}}$, it is easy to see that $C_{l+\frac {1}{2}}$ and $D_{l+\frac {1}{2}}$ are given by the
right formulae.
\end{proof}

\subsection[Case $k\in \frac{1}{2}+\mathbb{N}$]{Case $\boldsymbol{k\in \frac{1}{2}+\mathbb{N}}$}

\begin{proposition}\label{expl2}
If $\delta$ is non-critical and if $l$ is a~natural number, the symbol $S_{k-l}$ given in the proof of
Theorem~{\rm \ref{flatex}} is given by:
\begin{gather*}
\left[C'_l
\begin{pmatrix}
\partial_x^l&0
\\
0&\partial_x^l
\end{pmatrix}
+D'_l
\begin{pmatrix}
\partial_x^l&-\bar{D}_1\bar{D}_2\partial_x^{l-1}
\\
\bar{D}_1\bar{D}_2\partial_x^{l-1}&\partial_x^l
\end{pmatrix}
\right]S_k,
\end{gather*}
where the coefficients $C'_l$ and $D'_{l}$ are given in the following way:
\begin{gather*}
C'_l=\frac{\prod\limits_{i=0}^{l-1}{\Big(A_{k-\frac{1}{2}-i}B_{k-\frac{1}{2}-i}}\Big)}
{\prod\limits_{i=1}^l{(\alpha_k-\alpha_{k-i})}},
\qquad
D'_l=\frac{-l\prod\limits_{i=0}^{l-1}{\Big(A_{k-\frac{1}{2}-i}B_{k-\frac{1}{2}-i}\Big)}}
{2\Big(\alpha_k-\alpha_{k-\frac{1}{2}}\Big)\prod\limits_{i=1}^l{(\alpha_k-\alpha_{k-i})}}.
\end{gather*}

If $l\in\mathbb{N}+\frac {1}{2}$, we have:
\begin{gather*}
S_{k-l}= (-1)^{\tilde{S_k}+1}E'_l\left[
\begin{pmatrix}
B_{k-l}\bar{D_1}&B_{k-l}\bar{D}_2
\\
A_{k-l}\bar{D}_2&-A_{k-l}\bar{D}_1
\end{pmatrix}
\partial_x^{\lfloor l \rfloor}\right]S_k,
\end{gather*}
where $\lfloor l \rfloor$ is the integer part of $l$ and where the coefficients $E'_{l}$ are given
by the following formulae:
\begin{gather*}
E'_l=\frac{\prod\limits_{i=0}^{l-\frac{3}{2}}{\Big(A_{k-\frac{1}{2}-i}B_{k-\frac{1}{2}-i}\Big)}}
{2\Big(\alpha_k-\alpha_{k-\frac{1}{2}}\Big)\prod\limits_{i=1}^{l-\frac{1}{2}}{(\alpha_k-\alpha_{k-i})}}.
\end{gather*}

The first coefficient $E'_{\frac {1}{2}}$ is defined in the following way:
\begin{gather*}
E'_{\frac{1}{2}}=\frac{1}{2\left(\alpha_k-\alpha_{k-\frac{1}{2}}\right)}.
\end{gather*}
\end{proposition}
\begin{proof}
The proof is completely similar to the previous one.
\end{proof}

\begin{remark}
If $k\in\mathbb{N}$ \big(resp.\
{}$k\in\mathbb{N}+\frac{1}{2}$\big), if $l\in\mathbb{N}+\frac {1}{2}$ is between $\frac {3}{2}$ and $k-\frac
{1}{2}$ (resp.\
{}$k$) and if $\delta$ is such that the dif\/ference $\alpha_{k}-\alpha_{k-l}$ vanishes, the explicit
formula given in Proposition~\ref{expl1} (resp.\
Proposition~\ref{expl2}) is still well-def\/ined.
Since these formulae give $\mathfrak{spo}(2|2)$-equivariant quantizations and are continuous in
neighborhoods of these values of $\delta$, these formulae give also $\mathfrak{spo}(2|2)$-equivariant
quantizations for these values of $\delta$, by continuity.

If $k\in\mathbb{N}$ \big(resp.\
{}$k\in\mathbb{N}+\frac {1}{2}$\big), an $\mathfrak{spo}(2|2)$-equivariant quantization for symbols of degree
$k$ exists thus also for the values of $\delta$ equal to
\begin{gather*}
\frac {k^{2}-l^{2}+\frac {1}{4}}{2(k-l)}
\qquad
\left({\rm resp.\ }\;\frac {k^{2}-l^{2}-\frac {1}{4}}{2(k-l)}\right),
\end{gather*}
when the half-natural number $l$ (resp.\
the natural number $l$) varies from $\frac {1}{2}$ (resp.~$0$) to $k-\frac{3}{2}$.
\end{remark}

\begin{remark}
If $k\in\mathbb{N}$ and if $\alpha_{k}-\alpha_{k-\frac {1}{2}}=0$ (i.e.\
if $\delta=k$) or if $\alpha_{k}-\alpha_{k-1}=0$ \big(i.e.\
if $\delta=k-\frac {1}{2}$\big), then there is no $\mathfrak{spo}(2|2)$-equivariant quantization.
Indeed, in these situations, the equations whose the unknowns are $S_{k-\frac {1}{2}}$ and $S_{k-1}$ in the
system~\eqref{flatP} admit no solution because the left hand side of these equations vanishs whereas their
right hand side is not equal to zero.
Now, the existence of an $\mathfrak{spo}(2|2)$-equivariant quantization implies the existence of a~solution
of the system~\eqref{flatP}, as explained in the proof of Theorem~\ref{flatex}.
For the same reason, if $2\leqslant l\leqslant k$ and if $\alpha_{k}-\alpha_{k-l}=0$ \big(i.e.\
if $\delta=k-\frac {l}{2}$\big), then there is no $\mathfrak{spo}(2|2)$-equivariant quantization if the numbers
$B_{k-1},\ldots,B_{k-l+1}$ are dif\/ferent from zero \big(i.e.\
if $\lambda$ is dif\/ferent from the numbers $-\frac {k-1}{2},\ldots,-\frac {k-l+1}{2}$\big).
By cons, if one of these numbers vanishs, the system~\eqref{flatP} has inf\/initely many solutions.
In these situations, there are certainly inf\/initely many $\mathfrak{spo}(2|2)$-equivariant quantizations,
but our method linked to the Casimir operators does not allow to conclude.

Exactly in the same way, if $k\in\mathbb{N}+\frac {1}{2}$ and if $\alpha_{k}-\alpha_{k-\frac {1}{2}}=0$
\big(i.e.\
if $\delta=k-\frac {1}{2}$\big), there is no $\mathfrak{spo}(2|2)$-equivariant quantization.
If $1\leqslant l\leqslant k$ and if $\alpha_{k}-\alpha_{k-l}=0$ \big(i.e.\
if $\delta=k-\frac {l}{2}$\big), then there is no $\mathfrak{spo}(2|2)$-equivariant quantization if the numbers
$B_{k-\frac {1}{2}},\ldots,B_{k-l+\frac {1}{2}}$ are dif\/ferent from zero \Big(i.e.\
if $\lambda$ is dif\/ferent from the numbers $-\frac {k-\frac {1}{2}}{2},\ldots,-\frac {k-l+\frac{1}{2}}{2}$\Big).
As in the case $k\in\mathbb{N}$, if one of these numbers vanishs, there are surely inf\/initely many
$\mathfrak{spo}(2|2)$-equivariant quantizations, but the proof of their existence seems {\it a~priori}
dif\/f\/icult.
\end{remark}

\subsection*{Acknowledgements}

It is a~pleasure to thank T.~Leuther, P.~Mathonet, J.-P.~Michel and
V.~Ovsienko for numerous fruitful discussions and for their interest in our work.
We also warmly thank the referees for their suggestions and remarks which considerably improved the paper.
This research has been funded by the Interuniversity Attraction Poles Programme initiated by the Belgian
Science Policy Of\/f\/ice.

\pdfbookmark[1]{References}{ref}
\LastPageEnding


\begin{thebibliography}{99}
\footnotesize\itemsep=0pt

\bibitem{Ber87}
Berezin F.A., Introduction to superanalysis, \textit{Mathematical Physics and
  Applied Mathematics}, Vol.~9, D. Reidel Publishing Co., Dordrecht, 1987.

\bibitem{BHMP}
Boniver F., Hansoul S., Mathonet P., Poncin N., Equivariant symbol calculus for
  dif\/ferential operators acting on forms, \href{http://dx.doi.org/10.1023/A:1022251607566}{\textit{Lett. Math. Phys.}}
  \textbf{62} (2002), 219--232, \href{http://arxiv.org/abs/math.RT/0206213}{math.RT/0206213}.

\bibitem{BM}
Boniver F., Mathonet P., I{FFT}-equivariant quantizations, \href{http://dx.doi.org/10.1016/j.geomphys.2005.04.014}{\textit{J.~Geom.
  Phys.}} \textbf{56} (2006), 712--730, \href{http://arxiv.org/abs/math.RT/0109032}{math.RT/0109032}.

\bibitem{Bou00}
Bouarroudj S., Projectively equivariant quantization map, \href{http://dx.doi.org/10.1023/A:1007692910159}{\textit{Lett. Math.
  Phys.}} \textbf{51} (2000), 265--274, \href{http://arxiv.org/abs/math.DG/0003054}{math.DG/0003054}.

\bibitem{CapSil}
{\v{C}}ap A., {\v{S}}ilhan J., Equivariant quantizations for {AHS}-structures,
  \href{http://dx.doi.org/10.1016/j.aim.2010.01.016}{\textit{Adv. Math.}} \textbf{224} (2010), 1717--1734, \href{http://arxiv.org/abs/0904.3278}{arXiv:0904.3278}.

\bibitem{DLO}
Duval C., Lecomte P., Ovsienko V., Conformally equivariant quantization:
  existence and uniqueness, \textit{Ann. Inst. Fourier (Grenoble)} \textbf{49}
  (1999), 1999--2029, \href{http://arxiv.org/abs/math.DG/9902032}{math.DG/9902032}.

\bibitem{Fox}
Fox D.J.F., Projectively invariant star products, \href{http://dx.doi.org/10.1155/IMRP.2005.461}{\textit{Int. Math. Res. Pap.}}
   (2005), 461--510, \href{http://arxiv.org/abs/math.DG/0504596}{math.DG/0504596}.

\bibitem{GarMelOvs07}
Gargoubi H., Mellouli N., Ovsienko V., Dif\/ferential operators on supercircle:
  conformally equivariant quantization and symbol calculus, \href{http://dx.doi.org/10.1007/s11005-006-0129-8}{\textit{Lett. Math.
  Phys.}} \textbf{79} (2007), 51--65, \href{http://arxiv.org/abs/math-ph/0610059}{math-ph/0610059}.

\bibitem{Hansoul}
Hansoul S., Projectively equivariant quantization for dif\/ferential operators
  acting on forms, \href{http://dx.doi.org/10.1007/s11005-004-4293-4}{\textit{Lett. Math. Phys.}} \textbf{70} (2004), 141--153.

\bibitem{sarah}
Hansoul S., Existence of natural and projectively equivariant quantizations,
  \href{http://dx.doi.org/10.1016/j.aim.2007.03.007}{\textit{Adv. Math.}} \textbf{214} (2007), 832--864, \href{http://arxiv.org/abs/math.DG/0601518}{math.DG/0601518}.

\bibitem{KacAdvances}
Kac V.G., Lie superalgebras, \href{http://dx.doi.org/10.1016/0001-8708(77)90017-2}{\textit{Adv. Math.}} \textbf{26} (1977), 8--96.

\bibitem{Lecras}
Lecomte P.B.A., Classif\/ication projective des espaces d'op\'erateurs
  dif\/f\'erentiels agissant sur les densit\'es, \href{http://dx.doi.org/10.1016/S0764-4442(99)80211-0}{\textit{C.~R.~Acad. Sci. Paris
  S\'er.~I Math.}} \textbf{328} (1999), 287--290.

\bibitem{Leconj}
Lecomte P.B.A., Towards projectively equivariant quantization, \href{http://dx.doi.org/10.1143/PTPS.144.125}{\textit{Progr.
  Theoret. Phys. Suppl.}}  (2001), no.~144, 125--132.

\bibitem{LO}
Lecomte P.B.A., Ovsienko V.Yu., Projectively equivariant symbol calculus,
  \href{http://dx.doi.org/10.1023/A:1007662702470}{\textit{Lett. Math. Phys.}} \textbf{49} (1999), 173--196,
  \href{http://arxiv.org/abs/math.DG/9809061}{math.DG/9809061}.

\bibitem{LPS}
Leites D., Poletaeva E., Serganova V., On {E}instein equations on manifolds and
  supermanifolds, \href{http://dx.doi.org/10.2991/jnmp.2002.9.4.3}{\textit{J.~Nonlinear Math. Phys.}} \textbf{9} (2002),
  394--425, \href{http://arxiv.org/abs/math.DG/0306209}{math.DG/0306209}.

\bibitem{LMR}
Leuther T., Mathonet P., Radoux F., One
  {${\mathfrak{osp}}(p+1,q+1|2r)$}-equivariant quantizations, \href{http://dx.doi.org/10.1016/j.geomphys.2011.09.003}{\textit{J.~Geom.
  Phys.}} \textbf{62} (2012), 87--99, \href{http://arxiv.org/abs/1107.1387}{arXiv:1107.1387}.

\bibitem{LR}
Leuther T., Radoux F., Natural and projectively invariant quantizations on
  supermanifolds, \href{http://dx.doi.org/10.3842/SIGMA.2011.034}{\textit{SIGMA}} \textbf{7} (2011), 034, 12~pages,
  \href{http://arxiv.org/abs/1010.0516}{arXiv:1010.0516}.

\bibitem{MR}
Mathonet P., Radoux F., Natural and projectively equivariant quantizations by
  means of {C}artan connections, \href{http://dx.doi.org/10.1007/s11005-005-6783-4}{\textit{Lett. Math. Phys.}} \textbf{72} (2005),
  183--196, \href{http://arxiv.org/abs/math.DG/0606554}{math.DG/0606554}.

\bibitem{MR1}
Mathonet P., Radoux F., Cartan connections and natural and projectively
  equivariant quantizations, \href{http://dx.doi.org/10.1112/jlms/jdm030}{\textit{J.~Lond. Math. Soc.~(2)}} \textbf{76}
  (2007), 87--104, \href{http://arxiv.org/abs/math.DG/0606556}{math.DG/0606556}.

\bibitem{MR2}
Mathonet P., Radoux F., On natural and conformally equivariant quantizations,
  \href{http://dx.doi.org/10.1112/jlms/jdp024}{\textit{J.~Lond. Math. Soc.~(2)}} \textbf{80} (2009), 256--272,
  \href{http://arxiv.org/abs/0707.1412}{arXiv:0707.1412}.

\bibitem{MR3}
Mathonet P., Radoux F., Existence of natural and conformally invariant
  quantizations of arbitrary symbols, \href{http://dx.doi.org/10.1142/S1402925110001057}{\textit{J.~Nonlinear Math. Phys.}}
  \textbf{17} (2010), 539--556, \href{http://arxiv.org/abs/0811.3710}{arXiv:0811.3710}.

\bibitem{MR11}
Mathonet P., Radoux F., Projectively equivariant quantizations over the
  superspace {${\mathbb R}^{p|q}$}, \href{http://dx.doi.org/10.1007/s11005-011-0474-0}{\textit{Lett. Math. Phys.}} \textbf{98}
  (2011), 311--331, \href{http://arxiv.org/abs/1003.3320}{arXiv:1003.3320}.

\bibitem{Mel09}
Mellouli N., Second-order conformally equivariant quantization in dimension
  {$1|2$}, \href{http://dx.doi.org/10.3842/SIGMA.2009.111}{\textit{SIGMA}} \textbf{5} (2009), 111, 11~pages, \href{http://arxiv.org/abs/0912.5190}{arXiv:0912.5190}.

\bibitem{Mic09}
Michel J.-P., Quantif\/ication conform\'ement \'equivariante des f\/ibr\'es
  supercotangents, Ph.D. thesis, Universit\'e de la M\'editerran\'ee~-
  Aix-Marseille~II, 2009, available at
  \url{http://tel.archives-ouvertes.fr/tel-00425576}.

\bibitem{Mus97}
Musson I.M., On the center of the enveloping algebra of a classical simple
  {L}ie superalgebra, \href{http://dx.doi.org/10.1006/jabr.1996.7000}{\textit{J.~Algebra}} \textbf{193} (1997), 75--101.

\bibitem{Pin90}
Pinczon G., The enveloping algebra of the {L}ie superalgebra {${\rm
  osp}(1,2)$}, \href{http://dx.doi.org/10.1016/0021-8693(90)90265-P}{\textit{J.~Algebra}} \textbf{132} (1990), 219--242.

\bibitem{Ser99}
Sergeev A., The invariant polynomials on simple {L}ie superalgebras,
  \href{http://dx.doi.org/10.1090/S1088-4165-99-00077-1}{\textit{Represent. Theory}} \textbf{3} (1999), 250--280,
  \href{http://arxiv.org/abs/math.RT/9810111}{math.RT/9810111}.

\end{thebibliography}
\end{document}